\documentclass[english,draft|final]{siamltex}
\pdfoutput=1

\usepackage[latin9]{inputenc}
\usepackage{array}
\usepackage{verbatim}
\usepackage{float}
\usepackage{multirow}
\usepackage{amsmath}
\usepackage{amssymb}
\usepackage{graphicx}

\makeatletter

\providecommand{\tabularnewline}{\\}
\floatstyle{ruled}
\newfloat{algorithm}{tbp}{loa}
\providecommand{\algorithmname}{Algorithm}
\floatname{algorithm}{\protect\algorithmname}

\numberwithin{equation}{section}
\numberwithin{figure}{section}

\usepackage{algorithm,algpseudocode}
\usepackage{url}

\makeatother

\usepackage{babel}
\begin{document}

\title{Smoothed Hinge Loss and $\ell^{1}$ Support Vector Machines}

\author{Jeffrey Hajewski, Suely Oliveira and David E.~Stewart}

\date{$\today$}
\maketitle
\begin{abstract}
A new algorithm is presented for solving the soft-margin Support Vector
Machine (SVM) optimization problem with an $\ell^{1}$ penalty. This
algorithm is designed to require a modest number of passes over the
data, which is an important measure of its cost for very large data
sets. The algorithm uses smoothing for the hinge-loss function, and
an active set approach for the $\ell^{1}$ penalty. 
\end{abstract}
\begin{keywords}
support vector machine; smoothing; $\ell^{1}$ penalty
\end{keywords}
\begin{AMS}
Primary: 65K05, Secondary: 49M15
\end{AMS}

\section{Introduction}

\label{sec:Intro}\global\long\def\R{\mathbb{R}}
Dealing with large data sets has lead to a strong interest in methods
that have low iteration costs, such as stochastic gradient descent
(SGD) \cite{Niu:2011:HLA:2986459.2986537,robbins-munro1951}. More
classical methods such as Newton's method for optimization \cite[\S3.3]{nocedal}
are generally not used as their cost per iteration involves solving
linear systems which takes $\mathcal{O}(m^{3})$ operations where
$m$ is the number of unknowns. Contrary to conventional wisdom, we
argue that Newton's method and more sophisticated line search methods
are actually more appropriate for very large data problems, since
there the computational issues are typically due to the large number
of data items ($n$) rather than the dimension of the problem ($m$).
Wide data, where the dimension of the data vectors $\boldsymbol{x}_{i}$
is large compared to the number of data items, is still problematic
for Newton's method as the Hessian matrix is then singular. However,
in this paper we focus on $\ell^{1}$ Support Vector Machines ($\ell^{1}$SVMs)
and argue that Newton's method with a suitable line search and an
active-set strategy can also solve these problems very efficiently.
The algorithm developed here is, in part, inspired by \cite{MR1773265}
for the basis pursuit noise-reduction problem.

This algorithm and the numerical results are reported in the conference
paper \cite{hos:shll1svm}. The development of the line search algorithm
and justification for the convergence of the overall algorithm are
not reported in the conference paper.

The soft-margin SVM \cite[p.~263]{Duda:2000:PC:954544} for given
data $(\boldsymbol{x}_{i},y_{i})$, $i=1,\,2,\,\ldots,\,n$ where
each $y_{i}=\pm1$ minimizes
\begin{equation}
\frac{1}{2}\lambda\left\Vert \boldsymbol{w}\right\Vert _{2}^{2}+\frac{1}{n}\sum_{i=1}^{n}\max(0,\,1-y_{i}\,\boldsymbol{w}^{T}\boldsymbol{x}_{i})\label{eq:SVM-soft-margin}
\end{equation}
over all $\boldsymbol{w}\in\R^{m}$. Here the value of $\lambda>0$
is used to control the size of the vector $\boldsymbol{w}$. The function
$\max(0,\,1-y_{i}\,\boldsymbol{w}^{T}\boldsymbol{x}_{i})$ is called
the \emph{hinge-loss function} as it is based on the function $u\mapsto\max(0,\,u)$
whose graph looks like a hinge. The $\ell^{1}$SVM for the same data
minimizes
\begin{equation}
\frac{1}{2}\lambda\left\Vert \boldsymbol{w}\right\Vert _{2}^{2}+\frac{1}{n}\sum_{i=1}^{n}\max(0,\,1-y_{i}\,\boldsymbol{w}^{T}\boldsymbol{x}_{i})+\alpha\left\Vert \boldsymbol{w}\right\Vert _{1}\label{eq:L1-SVM}
\end{equation}
over $\boldsymbol{w}$. Here $\alpha>0$ controls the level of sparsity
of $\boldsymbol{w}$. Larger values tend to mean fewer components
of $\boldsymbol{w}$ are non-zero; if $\alpha$ is large enough then
$\boldsymbol{w}=0$. Note that this formulation is similar to, but
not the same as the 1-norm SVM of Zhu, Rosset, Hastie and Tibshirani
\cite{1norm-svm}. Also, the algorithm obtained here is $\mathcal{O}(n)$
with respect to the number of data points, while the algorithm of
Zhu et~al.\ is $\Omega(n^{2})$ as it involves identifying the intersections
of a descent line with the hyperplanes $1-y_{i}\,\boldsymbol{w}^{T}\boldsymbol{x}_{i}=0$
for each data point. Rather we use a smoothing approach for the sum
of the hinge-loss functions.

Traditionally, for optimization problems, the numbers of function,
gradient, and Hessian matrix evaluations are used to measure the cost
of the algorithm. For large-scale data mining types of optimization
problems, perhaps a different measure of performance is more important:
the number of passes over the data. The general form of most optimization
problems used in data mining is
\begin{equation}
\min_{\boldsymbol{w}}f(\boldsymbol{w}):=R(\boldsymbol{w})+\frac{1}{n}\sum_{i=1}^{n}\psi(\boldsymbol{x}_{i},y_{i};\boldsymbol{w})\label{eq:data-mining-gen-form}
\end{equation}
where $R$ is a regularization function, and $\psi$ is a loss function.
Provided $\boldsymbol{w}$ has relatively low dimension (say, below
$10^{3}$) and $n$ is large (say, $10^{5}$ to $10^{9}$), the cost
of computing $R(\boldsymbol{w})$ is modest and can be computed on
one processor, while the computations of $\psi(\boldsymbol{x}_{i},y_{i};\boldsymbol{w})$
should be carried out in parallel, and then summed via a parallel
reduction operation \cite{Blelloch:1996:PPA:227234.227246}. %

Computing the gradient of the objective function
\[
\nabla f(\boldsymbol{w})=\nabla R(\boldsymbol{w})+\frac{1}{n}\sum_{i=1}^{n}\nabla_{\boldsymbol{w}}\psi(\boldsymbol{x}_{i},y_{i};\boldsymbol{w}),
\]
which can be computed in a similar manner to the objective function,
except that the reduction (summation) is applied to the gradients
$\nabla_{\boldsymbol{w}}\psi(\boldsymbol{x}_{i},y_{i};\boldsymbol{w})$.
Similarly, the Hessian matrices can be computed in parallel but the
reduction (summation) is applied to the Hessian matrices $\text{Hess}_{\boldsymbol{w}}\psi(\boldsymbol{x}_{i},y_{i};\boldsymbol{w})$
of the loss functions. This may become an expensive step if $m$ becomes
large, as the reduction must be applied to objects of size $\mathcal{O}(m^{2})$
where $\boldsymbol{w}\in\R^{m}$. In such cases, a BFGS algorithm
may be appropriate instead of a direct Newton method.

If the function $R(\boldsymbol{w})$ is non-smooth in $\boldsymbol{w}$
(as is the case for (\ref{eq:L1-SVM})), then the optimization algorithm
needs to be adapted for it. 

\subsection{Problems with the line search.}

\label{subsec:Line-searches}Line searches are often needed in optimization
algorithms because the predicted step from Newton's method ``goes
too far'', or or in some other way results in an increase in the
objective function value or does not decrease it significantly. Suppose
the step for the Newton method is $\boldsymbol{d}$. If the quadratic
Taylor polynomial at $s=0$ to $f(\boldsymbol{w}+s\boldsymbol{d})$
is a poor approximation to $f(\boldsymbol{w}+s\boldsymbol{d})$, then
it may be necessary to perform many line search steps, which will
require many function evaluations. This is costly in the context of
parallel computation with high latency networks. 

Lack of smoothness can be a cause of this, and result in many costly
parallel reduction steps. Thus the shape of the $R(\boldsymbol{w})$
function must be known by the line search procedure at least to fairly
good accuracy. In the case of the $\ell^{1}$SVM problem, this means
that the nonsmoothness of the $\ell^{1}$ penalty must be explicitly
represented and used for the line search procedure. 

\subsection{Non-smoothness for $\ell^{1}$SVM.}

The advantage of using $\ell^{1}$SVM over a standard SVM formulation
is that the $\ell^{1}$ penalty tends to result in sparse solutions.
That is, with the $\ell^{1}$ penalty, the number of indexes $i$
where $w_{i}\neq0$ tends to be small. In fact, if the weight $\alpha>0$
is large enough, then the solution is $\boldsymbol{w}=0$. If $\alpha$
is smaller, we usually expect $w_{i}\neq0$ for a modest number of
indexes $i$. Sparse solutions have a number of advantages. There
is a much lower likelihood of over-fitting the data. The solution
is more likely to be ``explainable'' in the sense that the set of
$i$ where $w_{i}\neq0$ is small or modest, so that the method essentially
selects those parameters as being important. Finally, since fewer
parameters are used to create the ``fit'', there is probably less
noise in each of the parameters. Models with large numbers of parameters,
tend to have much less ``data per parameter'', so that the numerical
values obtained tend to be less reliable.

The disadvantage is that the numerical algorithm for performing the
optimization has to be adapted to deal with the non-smoothness. Since
the important non-smooth part of the objective function in (\ref{eq:L1-SVM})
is $\alpha\left\Vert \boldsymbol{w}\right\Vert _{1}$ is highly structured,
we can exploit this structure to create a fast and efficient algorithm.
To do this, an active set is maintained $\overline{\mathcal{I}}=\left\{ \,i\mid w_{i}\neq0\,\right\} $.
This needs to be expanded when new parameters $w_{i}$ are made active,
or available for optimization, and reduced when a line search indicates
that $w_{i}=0$ seems optimal for an active parameter $w_{i}$. If
$f(\boldsymbol{w})=g(\boldsymbol{w})+\alpha\left\Vert \boldsymbol{w}\right\Vert _{1}$
with $g$ smooth, an inactive parameter $w_{i}$ should be made active
if $\left|\partial g/\partial w_{i}(\boldsymbol{w})\right|>\alpha$.
With this strategy, many parameters can be made active in one step,
but only one active parameter can become inactive in one step. 

\subsection{Smoothing the hinge-loss function and convergence of Hessian matrices.}

\label{subsec:Smoothed-hinge-loss}
\begin{figure}
\begin{centering}
\includegraphics[width=0.8\textwidth]{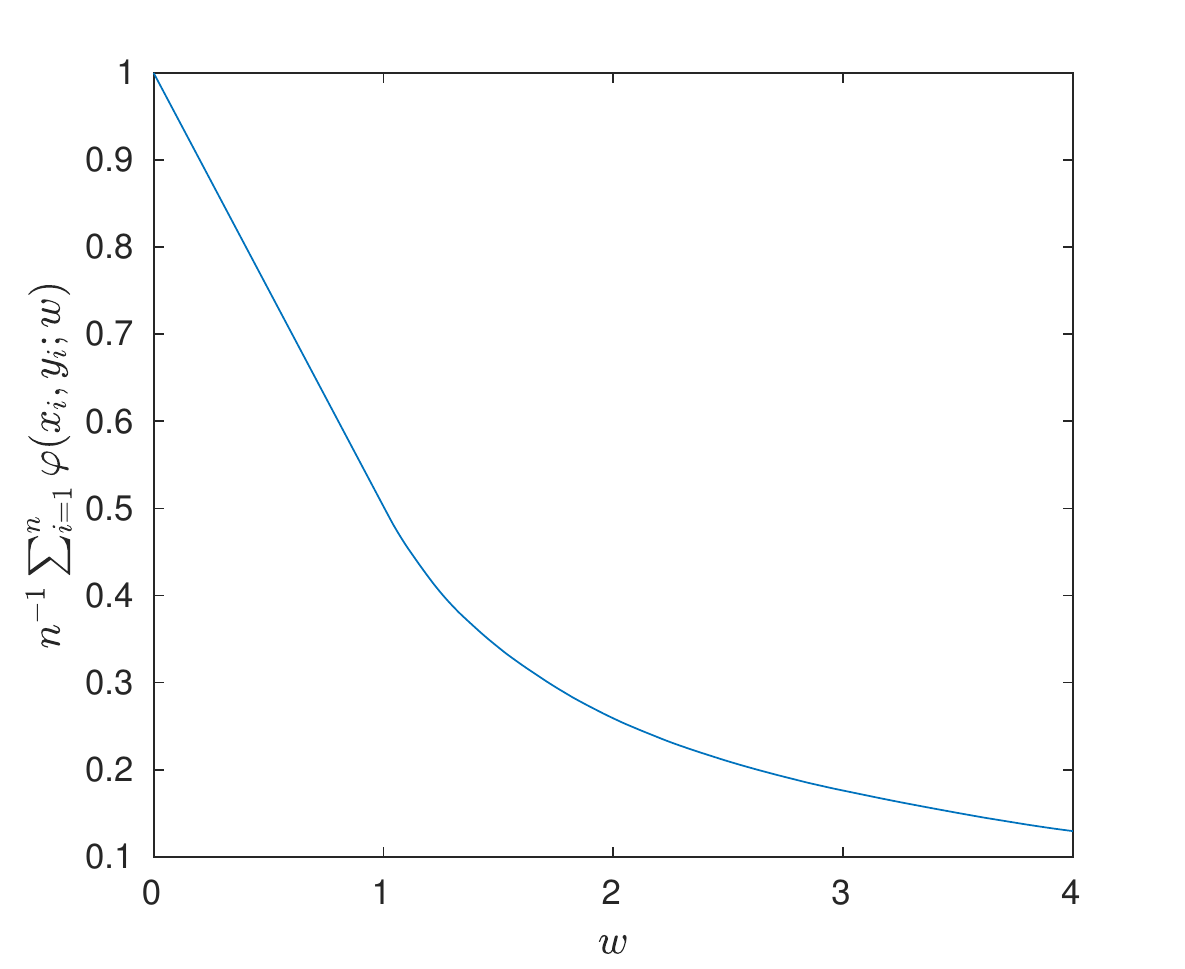}
\par\end{centering}
\caption{\label{fig:nearly-smooth}Plot of $n^{-1}\sum_{i=1}^{n}\psi(x_{i},y_{i};w)$
for randomly chosen data ($n=200$)}

\end{figure}

The hinge-loss function $\psi(\boldsymbol{x},y;\boldsymbol{w})=\max(0,\,1-y\,\boldsymbol{w}^{T}\boldsymbol{x})$
is a piecewise linear function of $\boldsymbol{w}$, and so its Hessian
matrix is zero or undefined. Thus
\[
\frac{1}{n}\sum_{i=1}^{n}\psi(\boldsymbol{x}_{i},y_{i};\boldsymbol{w})
\]
is also a piecewise linear function of $\boldsymbol{w}$, and thus
its Hessian matrix is either zero or undefined. On the other hand,
\[
\frac{1}{n}\sum_{i=1}^{n}\psi(\boldsymbol{x}_{i},y_{i};\boldsymbol{w})
\]
usually appears to be very smooth. For example, if $n=200$, and for
$y=+1$, $x$ is chosen randomly and uniformly from $[0,+1]$ while
for $y=-1$, $x$ is chosen uniformly and randomly from $[-1,0]$,
$n^{-1}\sum_{i=1}^{n}\psi(x_{i},y_{i};w)$ looks like Figure~\ref{fig:nearly-smooth}.

As $n\to\infty$ under some statistical assumptions detailed below,
the function $n^{-1}\sum_{i=1}^{n}\psi(x_{i},y_{i};w)$ approaches
a smooth function $h(\boldsymbol{w})$. Rather than compute the exact
Hessian matrix of $n^{-1}\sum_{i=1}^{n}\psi(x_{i},y_{i};w)$ with
respect to $\boldsymbol{w}$, which has very little to do with the
overall behavior of the function, we should aim to compute an approximation
to the Hessian matrix of $h(\boldsymbol{w})$. This can be done by
means of a smoothed hinge-loss function. Using a smoothed hinge-loss
function does not change the value of $n^{-1}\sum_{i=1}^{n}\psi(x_{i},y_{i};w)$
significantly, but does enable us to estimate the Hessian matrix of
$h(\boldsymbol{w})$, as well as its gradient.

For large data sets which come from some statistical distribution
with a $C^{2}$ probability density function, the mean of the hinge-loss
functions approaches a $C^{2}$ function
\begin{equation}
\frac{1}{n}\sum_{i=1}^{n}\psi(\boldsymbol{x}_{i},y_{i};\boldsymbol{w})\to\int\psi(\boldsymbol{x},+1;\boldsymbol{w})\,p_{1}(\boldsymbol{x})\,d\boldsymbol{x}+\int\psi(\boldsymbol{x},-1;\boldsymbol{w})\,p_{2}(\boldsymbol{x})\,d\boldsymbol{x}=:h(\boldsymbol{w})\label{eq:def-h(w)}
\end{equation}
as $n\to\infty$. Here $p_{1}(\boldsymbol{x})$ is the probability
density function of $\boldsymbol{x}$ given that $y=+1$, while $p_{2}(\boldsymbol{x})$
is the probability density function of $\boldsymbol{x}$ given that
$y=-1$. The values of the averages $(1/n)\sum_{i=1}^{n}\psi(\boldsymbol{x}_{i},y_{i};\boldsymbol{w})$
can then be well-approximated by a $C^{2}$ function. The difficulty
is in estimating the Hessian matrix of this unknown smooth function.
We can approximate $\psi(\boldsymbol{x},y;\boldsymbol{w})$ by a smoothed
hinge-loss function $\psi_{\epsilon}(\boldsymbol{x},y;\boldsymbol{w})$
given by
\begin{equation}
\psi_{\epsilon}(\boldsymbol{x},y;\boldsymbol{w})=\frac{1}{2}(u+\sqrt{\epsilon^{2}+u^{2}})\text{ where }u=1-y\,\boldsymbol{w}^{T}\boldsymbol{x}.\label{eq:def-smoothed-hinge-loss}
\end{equation}
What we want is that 
\[
\frac{1}{n}\sum_{i=1}^{n}\text{Hess}_{\boldsymbol{w}}\psi_{\epsilon}(\boldsymbol{x}_{i},y_{i};\boldsymbol{w})\approx\text{Hess}\,h(\boldsymbol{w})
\]
for $n$ sufficiently large. Now 
\[
\text{Hess}\,h(\boldsymbol{w})=\int\text{Hess}_{\boldsymbol{w}}\psi(\boldsymbol{x},+1;\boldsymbol{w})\,p_{1}(\boldsymbol{x})\,d\boldsymbol{x}+\int\text{Hess}_{\boldsymbol{w}}\psi(\boldsymbol{x},-1;\boldsymbol{w})\,p_{2}(\boldsymbol{x})\,d\boldsymbol{x}.
\]
and $\text{Hess}_{\boldsymbol{w}}\psi(\boldsymbol{x},y;\boldsymbol{w})=y^{2}\boldsymbol{x}\boldsymbol{x}^{T}\,\delta((1-y\,\boldsymbol{x}^{T}\boldsymbol{w})/\left\Vert y\boldsymbol{x}\right\Vert )$
where $\delta$ is the Dirac-$\delta$ distribution. That is,
\[
\int\text{Hess}_{\boldsymbol{w}}\psi(\boldsymbol{x},+1;\boldsymbol{w})\,p_{1}(\boldsymbol{x})\,d\boldsymbol{x}=\int_{\left\{ \boldsymbol{x}\mid1-\boldsymbol{w}^{T}\boldsymbol{x}=0\right\} }\boldsymbol{x}\boldsymbol{x}^{T}p_{1}(\boldsymbol{x})\,dS(\boldsymbol{x})
\]
where the latter is a surface integral. 

Note that $\psi(\boldsymbol{x},y;\boldsymbol{w})=j(1-y\boldsymbol{w}^{T}\boldsymbol{x})$
with $j(u)=\max(0,u)$, while $\psi_{\epsilon}(\boldsymbol{x},y;\boldsymbol{w})=j_{\epsilon}(1-y\boldsymbol{w}^{T}\boldsymbol{x})$
with $j_{\epsilon}(u)=\frac{1}{2}(u+\sqrt{\epsilon^{2}+u^{2}})$.
Now $j''(u)=\delta(u)$ while 
\begin{align*}
j_{\epsilon}''(u) & =\frac{1}{2}((\epsilon^{2}+u^{2})^{-1/2}-u^{2}(\epsilon^{2}+u^{2})^{-3/2})\\
 & =\frac{1}{2}\frac{\epsilon^{2}}{(\epsilon^{2}+u^{2})^{3/2}},
\end{align*}
which converges to $\delta(u)$ in the sense of distributions (and
the sense of measures, although weakly) as $\epsilon\to0$. Thus for
continuous $p_{1}$, 
\[
\int\text{Hess}_{\boldsymbol{w}}\psi_{\epsilon}(\boldsymbol{x},+1;\boldsymbol{w})\,p_{1}(\boldsymbol{x})\,d\boldsymbol{x}\to\int_{\left\{ \boldsymbol{x}\mid1-\boldsymbol{w}^{T}\boldsymbol{x}=0\right\} }\boldsymbol{x}\boldsymbol{x}^{T}p_{1}(\boldsymbol{x})\,dS(\boldsymbol{x})\qquad\text{as }\epsilon\to0.
\]
The integral on the right is an $(m-1)$-dimensional integral over
the hyperplane. Also, if we choose $\boldsymbol{x}_{i}$ independently,
and distributed according to the probability distribution $p_{1}$,
and the variance for the probability distribution $p_{1}$ is finite,
then by the Strong Law of Large Numbers \cite[p.~239]{loeve1977},
\begin{equation}
\frac{1}{n}\sum_{i=1}^{n}\text{Hess}_{\boldsymbol{w}}\psi_{\epsilon}(\boldsymbol{x}_{i},+1;\boldsymbol{w})\to\int\text{Hess}_{\boldsymbol{w}}\psi_{\epsilon}(\boldsymbol{x},+1;\boldsymbol{w})\,p_{1}(\boldsymbol{x})\,d\boldsymbol{x}\qquad\text{almost surely}.\label{eq:as-convergence-1}
\end{equation}
Since the same arguments apply for $p_{2}$ and the samples where
$y_{i}=-1$,
\begin{equation}
\lim_{\epsilon\to0}\lim_{n\to\infty}\frac{1}{n}\sum_{i=1}^{n}\text{Hess}_{\boldsymbol{w}}\psi_{\epsilon}(\boldsymbol{x}_{i},y_{i};\boldsymbol{w})=\text{Hess}_{\boldsymbol{w}}h(\boldsymbol{w})\qquad\text{almost surely}.\label{eq:as-convergence-2}
\end{equation}

To make this work in a practical sense, we need the number of samples
$n$ to be ``sufficiently large'' for a given $\epsilon>0$ in order
to have 
\[
\frac{1}{n}\sum_{i=1}^{n}\text{Hess}_{\boldsymbol{w}}\psi_{\epsilon}(\boldsymbol{x}_{i},y_{i};\boldsymbol{w})\approx\text{Hess}_{\boldsymbol{w}}h(\boldsymbol{w}),
\]
at least with high probability. A natural question is how large $n$
has to be for a given $\epsilon$ in order to have a good approximation.
Since $\psi_{\epsilon}(\boldsymbol{x},+1;\boldsymbol{w})$ only depends
on $\boldsymbol{w}^{T}\boldsymbol{x}$ (and similarly for $\psi_{\epsilon}(\boldsymbol{x},-1;\boldsymbol{w})$),
we only need $n\sim\text{const}\,\left\Vert \boldsymbol{w}\right\Vert \epsilon^{-1}$
as $\epsilon\to0$ in order to achieve a given level of accuracy in
approximating $\text{Hess}\,h(\boldsymbol{w})$. Thus the number of
data points needed to obtain a good approximation of the curvature
of the objective function is not exorbitant. 

\section{Development of the Algorithm}

\subsection{Choice of $\epsilon$}

We wish to use a modified Newton method for minimizing $f(\boldsymbol{w})$
from \ref{eq:data-mining-gen-form}. As noted in Section~\ref{sec:Intro},
$\text{Hess}\,f(\boldsymbol{w})$ is either undefined or $\text{Hess}\,R(\boldsymbol{w})$,
if it has a Hessian matrix. This is misleading, and will not lead
to fast convergence. Instead we use the smoothed hinge-loss function
$\psi_{\epsilon}$ for suitable $\epsilon>0$. The problem then is
to choose $\epsilon$. From the analysis in Section~\ref{subsec:Smoothed-hinge-loss},
we can choose $\epsilon$ to be inversely proportional to $n$, the
number of data points. With this approach, it might also be necessary
to adapt $\epsilon$ according to the distribution of the data points
$(\boldsymbol{x}_{i},y_{i})$. 

However, there is another approach which is agnostic regarding $n$
and the distribution of the data points. This is to simply begin with
a large value of $\epsilon$, minimize $f_{\epsilon}(\boldsymbol{w})$
over $\boldsymbol{w}$, then repeatedly reduce $\epsilon$ by (for
example) halving $\epsilon$, and then minimizing $f_{\epsilon}(\boldsymbol{w})$
over $\boldsymbol{w}$ with this new value of $\epsilon$. 

\subsection{Line-search algorithm}

In the context of parallel computing, it is important to keep the
number of function evaluations small. So it is important to use a
``good'' first guess. With Newton methods applied to smooth functions
$\psi$, it is traditional to use the step length $s=1$ with the
Newton step $\boldsymbol{d}=-(\text{Hess}\,\psi(\boldsymbol{w}))^{-1}\nabla\psi(\boldsymbol{w})$
followed by the Armijo line search (see \cite{arm:mflcfpd}, \cite[p.~33]{nw:no}).
However, with non-smooth functions, such as the $\ell^{1}$ penalty,
this choice can result in many function evaluations for a single line
search.

With the $\ell^{1}$ penalty, we have to consider the problem of minimizing
$\psi(\boldsymbol{w}+s\boldsymbol{d})+\alpha\left\Vert \boldsymbol{w}+s\boldsymbol{d}\right\Vert _{1}$
over $s\geq0$ efficiently where $\psi$ is a smooth function. Since
we can estimate the Hessian matrices accurately, we can use a quadratic
approximation for $\psi(\boldsymbol{w}+s\boldsymbol{d})\approx a\,s^{2}+b\,s+c$.
Then our line search seeks to minimize 
\begin{equation}
j(s):=a\,s^{2}+b\,s+c+\alpha\left\Vert \boldsymbol{w}+s\boldsymbol{d}\right\Vert _{1}\qquad\text{over }s\geq0.\label{eq:line-search-quad}
\end{equation}
Provided $a,\,\alpha\geq0$, this is a convex function, and so the
derivative $j'(s)$ is a non-decreasing function of $s$. Provided
$\left\Vert \boldsymbol{d}\right\Vert _{1}>b$ \emph{or} $a>0$ \emph{or}
$\alpha>0$, there is a global minimizer of $j$; if $a>0$ then it
is unique. The task is to compute this minimizer efficiently. This
minimizer is characterized by either $j'(s)=0$, or $j'(s^{-})\leq0$
and $j'(s^{+})\geq0$.

This can be done using a binary search algorithm that just uses the
data mentioned: $a$, $b$, $\alpha$, $\boldsymbol{w}$ and $\boldsymbol{d}$,
\emph{and no additional function evaluations}. All the information
needed from $\psi$ is $a$ and $b$, which can be computed from the
gradient and the Hessian matrix of $\psi$ at $\boldsymbol{w}$. 

Note that 
\begin{equation}
j'(s)=2as+b+\alpha\sum_{i=1}^{m}\text{sign}(w_{i}+sd_{i})\,d_{i}.\label{eq:j-prime}
\end{equation}
If $\alpha=0$ and $a>0$ then clearly the minimizing $s=-b/(2a)$.
Assuming $a>0$ and $\alpha\geq0$, the minimizing value of $s$ must
lie in the interval $[0,\,s_{max}]$ where $s_{max}=(\left|b\right|+\alpha\left\Vert \boldsymbol{d}\right\Vert _{1})/(2a)$. 

The points of discontinuity of $j'(s)$ are $\sigma_{i}=-w_{i}/d_{i}$,
$i=1,\,2,\,\ldots,\,m$. If any $d_{i}=0$, we can simply ignore $\sigma_{i}$.
Let $\left\{ \widehat{\sigma}_{1},\,\widehat{\sigma}_{2},\,\ldots,\,\widehat{\sigma}_{r}\right\} =\left\{ \sigma_{i}\mid\sigma_{i}>0\right\} $
with $\widehat{\sigma}_{1}<\widehat{\sigma}_{2}<\cdots<\widehat{\sigma}_{r}$.
Set $\widehat{\sigma}_{0}=0$. We check $j'(\widehat{\sigma}_{0}^{+})$
and $j'(\widehat{\sigma}_{r}^{+})$. If $j'(\widehat{\sigma}_{0}^{+})\geq0$
then the optimal $s$ is $s^{*}=0=\widehat{\sigma}_{0}$. If $a=0$
and $j'(\widehat{\sigma}_{r}^{+})<0$ then $j(s)\to-\infty$ as $s\to\infty$
and there is no minimum. If $a>0$ and $j'(\widehat{\sigma}_{r}^{+})<0$
then the optimal $s$ is
\[
s^{*}=-\frac{1}{2a}\left(b+\sum_{i=1}^{m}\text{sign}(d_{i})d_{i}\right)=-\frac{b+\alpha\left\Vert \boldsymbol{d}\right\Vert _{1}}{2a}=\widehat{\sigma}_{r}-\frac{j'(\widehat{\sigma}_{r}^{+})}{2a}>\widehat{\sigma}_{r},
\]
since $j'(\widehat{\sigma}_{r}^{+})=2a\widehat{\sigma}_{r}+b+\alpha\left\Vert \boldsymbol{d}\right\Vert _{1}$.%

Consider the sequence 
\begin{equation}
j'(\widehat{\sigma}_{0}^{+}),\,j'(\widehat{\sigma}_{1}^{-}),\,j'(\widehat{\sigma}_{1}^{+}),\,j'(\widehat{\sigma}_{2}^{-}),\,j'(\widehat{\sigma}_{2}^{+}),\,\ldots,\,j'(\widehat{\sigma}_{r}^{+}).\label{eq:j-prime-sequence}
\end{equation}
Since this sequence is a non-decreasing sequence, if $j'(\widehat{\sigma}_{0}^{+})<0$
and $j'(\widehat{\sigma}_{r}^{+})>0$ then it crosses from being $\leq0$
to $>0$ at some point. If $j'(\widehat{\sigma}_{i}^{\pm})=0$ for
some $i$ and choice of sign, then $s^{*}=\widehat{\sigma}_{i}$.
So we asume without loss of generality that $j'(\widehat{\sigma}_{i}^{\pm})\neq0$
for any $i$ and choice of sign. In this case, \emph{either} there
is an $i$ where $j'(\widehat{\sigma}_{i}^{-})<0$ and $j'(\widehat{\sigma}_{i}^{+})>0$,
\emph{or} there is an $i$ where $j'(\widehat{\sigma}_{i}^{+})<0$
and $j'(\widehat{\sigma}_{i+1}^{-})>0$. If $j'(\widehat{\sigma}_{i}^{-})<0$
and $j'(\widehat{\sigma}_{i}^{+})>0$, then $s^{*}=\widehat{\sigma}_{i}$.
If $j'(\widehat{\sigma}_{i}^{+})<0$ and $j'(\widehat{\sigma}_{i+1}^{-})>0$,
then $s^{*}\in(\widehat{\sigma}_{i},\,\widehat{\sigma}_{i+1})$. In
this latter case, for $s\in(\widehat{\sigma}_{i},\widehat{\sigma}_{i+1})$,
$j'(s)=j'(\widehat{\sigma}_{i}^{+})+2a(s-\widehat{\sigma}_{i})$,
so $s^{*}=\widehat{\sigma}_{i}-j'(\widehat{\sigma}_{i}^{+})/(2a)$. 

Finding the point where the sequence (\ref{eq:j-prime-sequence})
crosses zero can be carried out by binary search or a discrete version
of the bisection algorithm. Thus it can be computed in $\mathcal{O}(\log m)$
time as $r\leq m$. 

If the optimal value for $s^{*}$ is zero, then $\boldsymbol{d}$
is not a descent direction {[}Ref{]} and so some other direction should
be used. This can only occur if $\sigma_{i}=0$ for some $i$, indicating
that $w_{i}=0$. Then in this case, we need to remove $w_{i}$ from
the set of active variables.

\subsection{Combining the parts}

A complete algorithm is outlined in Algorithm~\ref{alg:SVM-l1-penalty-1-1}.
In this algorithm it should be noted that we use the following definitions:
\begin{align*}
\boldsymbol{a}\circ\boldsymbol{b} & =\boldsymbol{c}\qquad\text{where }c_{i}=a_{i}b_{i}\qquad(\text{Hadamard product})\\
\widehat{f}_{\alpha}(\boldsymbol{w}) & =(1/m)\sum_{i=1}^{m}\psi_{\alpha}(\boldsymbol{x}_{i},y_{i};1-y_{i}\,\boldsymbol{x}_{i}^{T}\boldsymbol{w})+\frac{1}{2}\lambda(\boldsymbol{w})^{T}\boldsymbol{w}\\
f_{\alpha}(\boldsymbol{w}) & =f_{\alpha}(\boldsymbol{w})+\mu\left\Vert \boldsymbol{w}\right\Vert _{1}.
\end{align*}
The inputs to $\psi_{\alpha}$ ($1-y_{i}\,\boldsymbol{x}_{i}^{T}\boldsymbol{w}$)
form the vector $\boldsymbol{e}-\boldsymbol{y}\circ(X\boldsymbol{w})$.
Note that $\boldsymbol{e}$ is the vector of 1's of the appropriate
size. This vector formulation is helpful in languages such as Matlab$^{TM}$.
Also, the matrix $X=[\boldsymbol{x}_{1},\,\boldsymbol{x}_{2},\,\ldots,\,\boldsymbol{x}_{n}]^{T}$
so that $X\boldsymbol{w}=[\boldsymbol{x}_{1}^{T}\boldsymbol{w},\,\ldots,\,\boldsymbol{x}_{n}^{T}\boldsymbol{w}]^{T}$.
Note that $\nabla\widehat{f}_{\alpha}(\boldsymbol{w})$ is well-defined
for all $\boldsymbol{w}$ provided $\alpha>0$, but that $f_{\alpha}$is
not smooth.

The algorithm used can be broken down into a number of pieces. At
the top level, the method can be considered as applying Newton's method
to a smoothed problem (smoothing parameter $\alpha$) keeping an inactive
set $\mathcal{I}=\left\{ \,i\mid w_{i}=0\,\right\} $. This inactive
set will need to change, either by gaining elements where $w_{j}\neq0$
but $w_{j}+sd_{j}=0$ resulting from the line search procedure, or
by losing elements where $w_{i}=0$ but the gradient component $g_{i}=\partial\widehat{f}_{\alpha}/\partial w_{i}(\boldsymbol{w})$
satisfies $\left|g_{i}\right|>\mu$ indicating that allowing $w_{i}\neq0$
will result in a lower objective function value. Note that $\widehat{f}_{\alpha}$
does not include the $\ell^{1}$ penalty term $\mu\left\Vert \boldsymbol{w}\right\Vert _{1}$.
The top-level computations are shown in Algorithm~\ref{alg:SVM-l1-penalty-1-1}.

An essential choice in this algorithm is \emph{not} to smooth the
$\ell^{1}$ penalty term, and instead use an active/inactive set approach.
If we had chosen to smooth the $\ell^{1}$ penalty term, then the
computational benefits of the smaller linear system in the Newton
step $\boldsymbol{d}_{\overline{\mathcal{I}}}\gets-H_{\overline{\mathcal{I}},\overline{\mathcal{I}}}^{-1}\widetilde{\boldsymbol{g}}_{\overline{\mathcal{I}}}$
would be lost. Instead, smoothing the $\ell^{1}$ term would mean
that the linear system to be solved would have size $m\times m$ where
$m$ is the dimension of $\boldsymbol{w}$. This would be particularly
important for problems with wide data sets where $m$ can be very
large. Instead, we expect that there would be bounds on the size of
$\overline{\mathcal{I}}$, the number of active weights $w_{i}\neq0$. 

\begin{algorithm}
\begin{algorithmic}[1]

\Require{$\alpha,\,\alpha_{min},\,\mu,\,\lambda>0$}

\Function{SVMsmooth}{$X,\,\boldsymbol{y},\,\boldsymbol{w},\,\lambda,\,\mu,\,\alpha,\,\alpha_{min}$}


\State$\mathcal{I}\gets\left\{ \,i\mid w_{i}=0\,\right\} $

\State$\mathcal{J}\gets\left\{ \,i\in\mathcal{I}\mid\left|g_{i}\right|>\mu\,\right\} $\Comment{Add
to active set}

\While{$\alpha>\alpha_{min}/\beta$}\Comment{While smoothing parameter
not at threshold}

  \State{}Carry out Newton step on smoothed problem

\EndWhile

\State$\mathbf{return}\;\boldsymbol{w}$

\EndFunction

\end{algorithmic}

\caption{\label{alg:SVM-l1-penalty-1-1}Algorithm for SVM with $\ell^{1}$
penalty}
\end{algorithm}

The Newton step computations are shown in Algorithm~\ref{alg:smoothed-l1-Newton-step}.
We first compute the gradient and the Hessian matrix. Care must be
taken at this point to ensure that we compute the correct gradient
for the components $j$ where $w_{j}=0$ but $\left|g_{j}\right|>\mu$.
The full Hessian matrix is not actually needed, just the ``active''
part of the Hessian matrix: $H_{\overline{\mathcal{I}},\overline{\mathcal{I}}}$.
The Newton step $\boldsymbol{d}$ is computed. If the predicted reduction
of the function value is sufficiently small, then we can assume the
problem for the current inactive set $\mathcal{I}$ and smoothing
parameter $\alpha>0$ has been solved to sufficient accuracy. Then
we can either reduce the current inactive set $\mathcal{I}$ or reduce
the smoothing parameter $\alpha$ as shown in Algorithm~\ref{alg:Adjust-active-set}.
The Newton steps then continue until either the inactive set or the
smoothing parameter is reduced. If the smoothing parameter goes below
$\alpha_{min}$, then the algorithm terminates.

\begin{algorithm}
\begin{algorithmic}[1]

  \State$\boldsymbol{g}\gets\nabla\widehat{f}_{\alpha}(\boldsymbol{w})$

  \State$\widetilde{\boldsymbol{g}}\gets\boldsymbol{g}+\mu\,\text{sign}(\boldsymbol{w})$

  \State$\widetilde{g}_{j}\gets\widetilde{g}_{j}+\mu\,\text{sign}(g_{j})$
for all $j\in\mathcal{J}$

  \State$H\gets\lambda\,I+(1/m)X^{T}\text{diag}(\psi_{\alpha}''(\boldsymbol{z}))X$\Comment{$H$
is Hessian matrix}

  \State$\boldsymbol{d}_{\overline{\mathcal{I}}}\gets-H_{\overline{\mathcal{I}},\overline{\mathcal{I}}}^{-1}\widetilde{\boldsymbol{g}}_{\overline{\mathcal{I}}}$;
$\boldsymbol{d}_{\mathcal{I}}\gets0$\Comment{Newton step}

  \If{$\left|\boldsymbol{d}^{T}\widetilde{\boldsymbol{g}}\right|<\alpha/10$}\Comment{If
smoothed problem nearly solved for $\alpha$ and $\mathcal{I}$\ldots{}}

    \State{}Adjust active set \& reduce smoothing parameter

    \State$\mathbf{continue}$

  \EndIf

  \State$s\gets\text{LinesearchL1}(\boldsymbol{w},\boldsymbol{d},\boldsymbol{g}^{T}\boldsymbol{d},\frac{1}{2}\boldsymbol{d}^{T}H\boldsymbol{d},\mu)$;
$\boldsymbol{w}^{+}\gets\boldsymbol{w}+s\,\boldsymbol{d}$

  \While{$f_{\alpha}(\boldsymbol{w}^{+})>f_{\alpha}(\boldsymbol{w})+c_{1}\,s\,\boldsymbol{d}^{T}\widetilde{\boldsymbol{g}}$}\Comment{Armijo
line search}

    \State$s\gets s/2$; $\boldsymbol{w}^{+}\gets\boldsymbol{w}+s\,\boldsymbol{d}$

  \EndWhile

  \State$\boldsymbol{w}\gets\boldsymbol{w}^{+}$; $\mathcal{I}\gets\left\{ \,i\mid w_{i}=0\,\right\} $\Comment{Add
to $\mathcal{I}$ if line search indicates}

\end{algorithmic}

\caption{\label{alg:smoothed-l1-Newton-step}Newton step}
\end{algorithm}

\begin{algorithm}
\begin{algorithmic}[1]

    \State$\mathcal{J}'\gets\left\{ \,i\in\mathcal{I}\mid\left|g_{i}\right|>\mu\,\right\} $

    \If{$\mathcal{J}'\neq\mathcal{J}$}

      \State$\mathcal{J}\gets\mathcal{J}'$; $\mathcal{I}\gets\mathcal{I}\backslash\mathcal{J}'$;
$\mathbf{continue}$

    \EndIf

    \State$\alpha\gets\alpha/\beta$\Comment{Reduce $\alpha$ and
optimize for this new $\alpha$}

\end{algorithmic}

\caption{\label{alg:Adjust-active-set}Adjust active set \& reduce smoothing
parameter}
\end{algorithm}

The linesearch algorithm is shown as Algorithm~\ref{alg:l1-linesearch}.

\begin{algorithm}
\begin{algorithmic}[1]

\Require{$a,\,\mu\geq0$ and $\boldsymbol{d}\neq0$ and either $a>0$
or $\mu\left\Vert \boldsymbol{d}\right\Vert _{1}>-b$}

\Require{$s_{max}\geq0$}

\Function{LinesearchL1}{$\boldsymbol{w},\boldsymbol{d},b,a,\mu,s_{max}$}

\Comment{returns $s$ that minimizes $a\,s^{2}+b\,s+\mu\left\Vert \boldsymbol{w}+s\boldsymbol{d}\right\Vert _{1}$
over $0\leq s\leq s_{max}$}

  \State$n\gets\text{dimension}(\boldsymbol{w})$

  \State{}find function $p\colon\left\{ 1,2,\ldots,m\right\} \to\left\{ 1,2,\ldots,n\right\} $
where 

  \State$\qquad$$\text{range}(p)=\left\{ \,j\mid-w_{j}/d_{j}>0\,\right\} $
and $[-w_{p(i)}/d_{p(i)}]_{i=1}^{m}$ is sorted

  \State$i_{1}\gets0$; $j_{1}\gets0$; $s_{1}\gets0$; $\mathit{slope}_{1}=b+\mu\,\text{sign}(\boldsymbol{w})^{T}\boldsymbol{d}$

  \State$i_{2}\gets m+1$; $j_{2}\gets n+1$; $s_{2}\gets+\infty$;
$\mathit{slope}_{2}\gets\begin{cases}
+\infty, & \text{if }a>0,\\
b+\mu\left\Vert \boldsymbol{d}\right\Vert _{1}, & \text{if }a=0.
\end{cases}$

  \State\textbf{if }$\mathit{slope}_{1}\geq0$ \textbf{then} $\mathbf{return}\;s_{1}$
\textbf{end if}

  \State\textbf{if }$\mathit{slope}_{2}\leq0$ \textbf{then} $\mathbf{return}\;s_{2}$
\textbf{end if}

\While{$i_{2}>i_{1}+1$}\Comment{binary search}

  \State$i\gets\left\lfloor (i_{1}+i_{2})/2\right\rfloor $; $j\gets p(i)$

  \State$s\gets-w_{j}/d_{j}$\Comment{compute slopes on either side
of $s$}

  \State$\mathit{slope}_{0}\gets2as+b+\sum_{k:k\neq j}\text{sign}(w_{k}+sd_{k})d_{k}$

  \State$\mathit{slope}_{+}\gets\mathit{slope}_{0}+\mu\left|d_{j}\right|$;
$\mathit{slope}_{-}\gets\mathit{slope}_{0}-\mu\left|d_{j}\right|$

  \If{($\mathit{slope}_{-}=0$ or $\mathit{slope}_{+}=0$) or ($\mathit{slope}_{-}<0$
and $\mathit{slope}_{+}>0$)}

    \State$\mathbf{return}\;s$

  \Else{} \textbf{if} $\mathit{slope}_{+}<0$\textbf{ then} $i_{1}\gets i$
\textbf{else} $i_{2}\gets i$ \textbf{end if}

  \EndIf

\EndWhile

\Comment{Note that $i_{2}=i_{1}+1$ \& the optimal $s$ is in $(s_{1},s_{2})$}

\State$s\gets(s_{1}\mathit{slope}_{2}-s_{2}\mathit{slope}_{1})/(\mathit{slope}_{2}-\mathit{slope}_{1})$

\EndFunction

\end{algorithmic}

\caption{\label{alg:l1-linesearch}Linesearch algorithm for quadratic plus
$\ell^{1}$ penalty}
\end{algorithm}

The actual implementation differs slightly from the pseudo-code in
that the recomputation of $\mathcal{I}$ on line~?? of Algorithm~\ref{alg:smoothed-l1-Newton-step}
uses some additional information returned from \textsc{LinesearchL1}:
in floating point arithmetic there is no guarantee that $\mathcal{I}\gets\left\{ \,i\mid w_{i}=0\,\right\} $
will identify components $w_{i}$ that would be set to zero in exact
arithmetic. Specifically, setting $s\gets-w_{j}/d_{j}$ in does not
ensure that $w_{j}+s\,d_{j}$ evaluates to zero in floating point
arithmetic. So the linesearch function \textsc{LinesearchL1} actually
returns both $s$ and $j_{1}$ and $j_{2}$: if $j_{1}=j_{2}$, then
$s=-w_{j}/d_{j}$ for $j=j_{1}=j_{2}$ and we would set $w_{j}+s\,d_{j}=0$
and the new set $\mathcal{I}$ is the old $\mathcal{I}$ plus~$j$. 

Thus, elements can be added to $\mathcal{I}$ (line~14 of Algorithm~\ref{alg:smoothed-l1-Newton-step})
as well as removed from $\mathcal{I}$ (line~3 of Algorithm~\ref{alg:Adjust-active-set}).
Note, however, that while this approach can remove multiple elements
of $\mathcal{I}$ in a single iteration, only a single element can
be added per iteration. This means that the dimension of $\boldsymbol{w}$
can strongly affect the number of iterations if $\mathcal{I}$ at
the optimum has many elements. As removal of elements of $\mathcal{I}$
is easier than addition of elements, it is probably better to begin
with $\mathcal{I}=\left\{ 1,2,\ldots,m\right\} $ and $\boldsymbol{w}=0$. 

\section{Results\label{sec:results}}

Our experiments explore SmSVM's performance using both real and synthetic
data (see Table~\ref{tab:data} for a detailed description of the
data). We look at the ability of our models to accurately classify
test data while maintaining, and in many cases improving, state of
the art training time. Additionally, we study the robustness of the
model as the training data becomes increasingly sparse by increasing
the number of components equal to zero in the two centroids used to
generate the synthetic data. This is discussed in greater detail in
Section \ref{section:data}. The results of this Section were previously
published in \cite{hos:shll1svm}.

\begin{table}[h!]
\caption{Description of datasets used in performance comparison experiments.
Sparsity refers to the percentage of the data with a value of 0. Note
that for the synthetic datasets, the sparsity levels vary based on
experiment.}
\label{tab:data}
\centering{}%
\begin{tabular}{lrrc}
\hline 
\textbf{Name}  & \textbf{Count}  & \textbf{Dimension}  & \textbf{Sparsity}\tabularnewline
\hline 
Australian  & 690  & 14  & 13\%\tabularnewline
Colon Cancer  & 62  & 2,000  & 0\%\tabularnewline
CoverType  & 581,012  & 54  & 78\%\tabularnewline
Synthetic (tall)  & 10,000  & 50  & N/A\tabularnewline
Synthetic (wide)  & 50  & 2,500  & N/A\tabularnewline
\hline 
\end{tabular}
\end{table}

We compare our algorithms against conjugate gradient (Polak-Ribière
Plus~\cite{nocedal,polakribiere}), subgradient descent, stochastic
subgradient descent, and coordinate descent (via LIBLINEAR~\cite{liblinear08}).
In the case of conjugate gradient, since our loss function is non-smooth,
we use a subgradient in place of the gradient, where a subgradient
is any element of the subdifferential~\cite{baptiste2004}: 
\[
\partial f(\boldsymbol{x})=\left\{ \,\boldsymbol{g}\in\mathbb{R}^{n}\mid f(\boldsymbol{y})\geq f(\boldsymbol{x})+\boldsymbol{g}^{T}(\boldsymbol{y}-\boldsymbol{x})\ \ \forall\boldsymbol{y}\in\mathbb{R}^{n}\,\right\} .
\]
Table~\ref{tab:names} describes the naming convention used in the
following sections along with a brief description of the algorithms.

\begin{table}[h!]
\centering{}\caption{Summary of objectives and algorithms.}
\label{tab:names} \vspace{1em}
\begin{tabular}{lp{0.6\linewidth}}
\hline 
\textbf{Name } & \textbf{Description}\tabularnewline
\hline 
SmSVM\textendash $\ell^{2}$  & $\ell^{2}$ regularization\tabularnewline
SmSVM\textendash $\ell^{1}$\textendash $\ell^{2}$  & $\ell^{2}$ and $\ell^{1}$ regularization\tabularnewline
LinearSVC  & LIBLINEAR~\cite{liblinear08}\tabularnewline
SGD $\ell^{2}$  & SGD $\ell^{2}$ regularization\tabularnewline
SSGD $\ell^{2}$ mb  & SGD $\ell^{2}$ regularization mini-batch size of 32\tabularnewline
CG  & Polak-Ribière Plus~\cite{polakribiere} conjugate gradient solves
equation (\ref{eq:cgplain})\tabularnewline
CG \textendash{} $\ell^{2}$  & Polak-Ribière Plus~\cite{polakribiere} conjugate gradient with $\ell^{2}$
regularization\tabularnewline
\hline 
\end{tabular}
\end{table}

We consider four different optimization problems in the following
experiments. SmSVM\textendash $\ell^{2}$ and CG $\ell^{2}$ solve
the optimization problem defined by equation (\ref{eq:SVM-soft-margin}),
while SmSVM\textendash $\ell^{1}$\textendash $\ell^{2}$ minimizes
the loss function defined in equation (\ref{eq:L1-SVM}). The standard
conjugate gradient optimizer minimizes (\ref{eq:cgplain}). 
\begin{equation}
\frac{1}{n}\sum_{i=1}^{n}\max\{0,1-y_{i}\boldsymbol{w}^{T}\boldsymbol{x}_{i}\}\label{eq:cgplain}
\end{equation}

SGD $\ell^{2}$ minimizes (\ref{eq:SVM-soft-margin}) using a stochastic
gradient descent method \cite{1606.04838,robbins-munro1951}. The
LinearSVC model, which is a Python wrapper over LIBLINEAR provided
by Scikit-learn \cite{scikit-learn}, solves a scaled version of equation
(\ref{eq:SVM-soft-margin}), shown in equation (\ref{eq:liblinear}).
\begin{equation}
C\sum_{i=1}^{n}\max\{0,1-y_{i}\boldsymbol{w}^{T}\boldsymbol{x}_{i}\}+\frac{1}{2}\left\Vert \boldsymbol{w}\right\Vert _{2}^{2},\hspace{0.2cm}C>0\label{eq:liblinear}
\end{equation}
In our case, this is optimized via coordinate descent (see \cite{liblinear08}
for details).

For experiments involving synthetic data, new data is generated each
repetition of the experiment. Unless otherwise noted, all experiments
are performed 50 times.


\subsection{Data}

\label{section:data} We use both synthetic and real data to compare
the SmSVM algorithms against the conjugate gradient and gradient descent
algorithms mentioned in Table~\ref{tab:names}. Table~\ref{tab:data}
describes the data used in the experiments. The synthetic data is
generated by creating two centroids with components randomly sampled
from $N(0,1)$, scaling the centroids, and then sampling $\boldsymbol{x}\sim N(\boldsymbol{c}_{i},\mathbb{I}_{m})$
where $\boldsymbol{c}_{i}\in\mathbb{R}^{m}$ is the respective centroid.
Sparse data is created by setting randomly selected components of
the centroids to zero, and then randomly sampling about the updated
centroids. The real datasets used in the experiments were sourced
from the UCI Machine Learning Repository \cite{Dua:2017}. The Australian
and Colon Cancer \cite{coloncancer99} datasets were chosen for their
shapes, with the Australian dataset being tall and narrow while the
Colon Cancer dataset is short and wide. The CoverType \cite{covtype}
dataset was chosen due to its size and is the largest dataset we ran
in our experiments. As noted in Table~\ref{tab:results}, the CoverType
dataset was only run 20 times, due to compute time constraints.

Table~\ref{tab:results} summarizes the overall results of test accuracy
and training time on the four datasets.

\begin{table*}
\centering{}\caption{Numerical results for the real world datasets. These results are the
average of 50 independent runs.}
\label{tab:results} %
\noindent\begin{minipage}[t]{1\columnwidth}%
\begin{center}
\begin{tabular}{l|cc|cc|cc}
\hline 
\multirow{2}{*}{Algorithm} & \multicolumn{2}{c|}{Australian} & \multicolumn{2}{c|}{Colon Cancer} & \multicolumn{2}{c}{CoverType \footnote{Results based on 20 runs due to computational requirements.}}\tabularnewline
 & Acc.  & Time (s)  & Acc.  & Time (s)  & Acc.  & Time (s) \tabularnewline
\hline 
SmSVM  & 44.5  & 0.051  & 38.0  & 5.143  & 51.2  & 44.2 \tabularnewline
SmSVM\textendash $\ell^{2}$  & 85.9  & 0.058  & 66.3  & 32.851  & 69.8  & 148.7 \tabularnewline
SmSVM\textendash $\ell^{1}$\textendash $\ell^{2}$  & \textbf{86.1 } & \textbf{0.002 } & \textbf{84.0 } & 0.918  & 69.5  & \textbf{1.1 }\tabularnewline
LinearSVC-Hinge  & 85.2  & 0.007  & 66.9  & \textbf{0.008 } & \textbf{76.3 } & 182.9 \tabularnewline
SGD $\ell^{2}$  & 85.9  & 0.008  & \textbf{84.0 } & 0.023  & 68.3  & 30.7 \tabularnewline
SSGD $\ell^{2}$ mb  & 85.9  & 0.058  & 80.9  & 0.016  & 63.9  & 53.0 \tabularnewline
CG  & 86.0  & 1.375  & 75.1  & 0.940  & 68.4  & 754.6 \tabularnewline
CG \textendash{} $\ell^{2}$  & 85.9  & 1.355  & 77.1  & 2.350  & 68.4  & 746.3 \tabularnewline
\hline 
\end{tabular}
\par\end{center}
\begin{center}
\begin{tabular}{l|cc|cc}
\hline 
\multirow{2}{*}{Algorithm} & \multicolumn{2}{c|}{Synthetic Tall} & \multicolumn{2}{c}{Synthetic Wide}\tabularnewline
 & Acc.  & Time (s)  & Acc.  & Time (s)\tabularnewline
\hline 
SmSVM  & 49.9  & 0.15  & \textbf{100 } & 11.89\tabularnewline
SmSVM\textendash $\ell^{2}$  & 84.3  & 0.16  & \textbf{100 } & 23.41\tabularnewline
SmSVM\textendash $\ell^{1}$\textendash $\ell^{2}$  & 76.0  & \textbf{0.01 } & 93.6  & 0.21\tabularnewline
LinearSVC-Hinge  & \textbf{100 } & 0.03  & \textbf{100 } & \textbf{0.01}\tabularnewline
SGD $\ell^{2}$  & 61.1  & 0.16  & 51.2  & 0.02\tabularnewline
SSGD $\ell^{2}$ mb  & 77.5  & 0.62  & 54.0  & 0.02\tabularnewline
CG  & 94.9  & 4.84  & 82.0  & 0.19\tabularnewline
CG \textendash{} $\ell^{2}$  & 78.6  & 5.64  & 80.8  & 1.29\tabularnewline
\hline 
\end{tabular}
\par\end{center}%
\end{minipage}
\end{table*}


\subsection{Performance and Implementation}

As seen in Table~\ref{tab:results}, SmSVM\textendash $\ell^{1}$\textendash $\ell^{2}$
performs well across a variety of dataset types, and is beaten only
by other SmSVM algorithms and LIBLINEAR \cite{liblinear08}. Most
notable is the incredibly fast training time, which is due to the
optimizations made available via the feature selection property of
the $\ell^{1}$ norm. We optimize the matrix-vector and vector-vector
operations by reducing the problem size to that of the active set
dimension. The reduction in problem size yields substantial computational
savings in problems where the active-set is small. The savings are
apparent in the real-world datasets where SmSVM\textendash $\ell^{1}$\textendash $\ell^{2}$
finished training, in some cases, by an order of magnitude shorter
time. One interesting aspect of SmSVM\textendash $\ell^{1}$\textendash $\ell^{2}$'s
performance is its apparent struggle in terms of training time on
the Colon Cancer dataset, which is a dense dataset. Although SmSVM\textendash $\ell^{1}$\textendash $\ell^{2}$
tied with SGD $\ell^{2}$ for the top test accuracy, the SmSVM family
of algorithms were among the slowest to finish training.

Perhaps the most surprising result is the performance on the CoverType
\cite{covtype} dataset. Consisting of nearly 600,000 data points
and roughly 70MB in uncompressed libSVM sparse format (only non-zero
values and their indices are given, everything else is assumed 0).
LIBLINEAR took nearly 3 minutes to train on this dataset, achieving
a best-in-class test accuracy, while SmSVM\textendash $\ell^{1}$\textendash $\ell^{2}$
trained in just over one second and achieving nearly a second-place
test accuracy. The closest algorithm to SmSVM\textendash $\ell^{1}$\textendash $\ell^{2}$
in terms of training time is SGD $\ell^{2}$, which was nearly 30
seconds slower and had a lower test accuracy.

The SmSVM, CG, and SGD optimizers were all implemented in pure python
and make extensive use of Numpy \cite{numpy}. We implemented these
algorithms as efficiently as possible, and in particular, focused
on reducing data-copying as much as possible. The LIBLINEAR implementation
was accessed via Scikit-learn \cite{scikit-learn}, which provides
a Python wrapper on the C++ implementation. %
{} 

\section{Discussion}

We have introduced SmSVM, a new approach to solving soft-margin SVM,
which is capable of strong test accuracy without sacrificing training
speed. This is achieved by smoothing the hinge-loss function and using
an active set approach to the the $\ell^{1}$ penalty. SmSVM provides
improved test accuracy over LIBLINEAR with comparable, and in some
cases reduced, training time. SmSVM uses orders of magnitude fewer
gradient calculations and a modest number of passes over the data
to achieve its results, meaning it will scales well for increasing
problem sizes. SmSVM\textendash $\ell^{1}$\textendash $\ell^{2}$
optimizes its matrix-vector and vector-vector calculations by reducing
the problem size to that of the active set. For even modestly sized
problems this results in significant savings with respect to computational
complexity.

Overall the results are quite promising. On the real and synthetic
datasets, our algorithms outperform or tie the competition in test
accuracy 80\% of the time and have the fastest training time 60\%
of the time. The time savings are increasingly significant as the
number of data points grows. The results of the wide synthetic dataset
are somewhat surprising in that the SmSVM\textendash $\ell^{1}$\textendash $\ell^{2}$
algorithm performed worse than the SmSVM\textendash $\ell^{2}$ algorithm
with respect to test accuracy. This is likely due to the SmSVM\textendash $\ell^{1}$\textendash $\ell^{2}$
algorithm pushing features out of the active set too aggressively.
On the other hand, training time was nearly two orders of magnitude
faster, due to the active set being considerably smaller, which allows
us to optimize some of the linear algebra operations.

SmSVM is implemented in Python, making it easy to modify and understand.
The use of Numpy keeps linear algebra operations optimized\textendash this
is important when competing against frameworks such as LIBLINEAR,
which is implemented in C++. Testing SmSVM on larger datasets, incorporating
GPU acceleration to the linear algebra, and exploring distributed
implementations are promising future directions. 


\def\cprime{$'$} \def\cprime{$'$}

\end{document}